\documentclass[12pt]{amsart}
\usepackage[foot]{amsaddr}  

\usepackage{enumerate}
\usepackage{amsmath, amssymb, amsthm, amsfonts}
\usepackage{tikz}
\usepackage{geometry}
\usepackage[pdftex, bookmarks=true]{hyperref}
\usepackage{bbm} %for the 1 symbol of indicator functions
\usepackage{url}
\usepackage{cmap} %supports pdf copy-paste
\usepackage{verbatim}  %for comment environment
\usepackage{bm}  %for bold in math mode

\usepackage{marginnote}

\newtheorem{introthm}{Theorem}

\newtheorem{theorem}{Theorem}[section]
\newtheorem{proposition}[theorem]{Proposition}
\newtheorem{corollary}[theorem]{Corollary}
\newtheorem{lemma}[theorem]{Lemma}
\newtheorem{conjecture}[theorem]{Conjecture}

\newtheorem*{claim*}{Claim}

\theoremstyle{definition}
\newtheorem{remark}[theorem]{Remark}
\newtheorem*{remark*}{Remark}

\newtheorem{definition}[theorem]{Definition}
\newtheorem{question}[theorem]{Question}

\newtheorem{example}[theorem]{Example}

\def\ds{\displaystyle}

\def\E{\mathbb{E}}  %expectation
\def\P{\mathbb{P}}  %probability

\def\Rb{\mathbb{R}}
\def\Zb{\mathbb{Z}}
\def\Nb{\mathbb{N}}
\def\Ac{\mathcal{A}}
\def\Bc{\mathcal{B}}

\def\Ic{\mathcal{I}}

\def\Pc{\mathcal{P}}  %power set

\def\al{\alpha}

\def\de{\delta}
\def\eps{\varepsilon}
\def\la{\lambda}

\def\sm{\setminus}

\DeclareMathOperator{\VP}{VP}
\DeclareMathOperator{\VPcnt}{VP_{cnt}}
\DeclareMathOperator{\VPfull}{VP_{full}}
\DeclareMathOperator{\VPcl}{\overline{VP}}

\DeclareMathOperator{\conv}{Conv} %convex hull
\DeclareMathOperator{\diam}{diam}

\newcommand{\defeq}{\mathrel{\vcenter{\baselineskip0.5ex \lineskiplimit0pt
                     \hbox{\scriptsize.}\hbox{\scriptsize.}}}%
                     =}

\def\ind{\mathbbm{1}} %the 1 symbol of indicator functions

\def\dd{\, \mathrm{d}}
\def\fpi{\varphi_\pi}
\def\Hupp{\overline{H}}
\def\Hlow{\underline{H}}
\def\Hfin{\overline{H}_{\mathrm{fin}}}
\def\Hcnt{\overline{H}_{\mathrm{cnt}}}
\def\Hfull{\overline{H}_{\mathrm{full}}}
\def\chif{\chi_{\mathrm{frac}}}
\def\omf{\omega_{\mathrm{frac}}}

\DeclareFontFamily{U}{matha}{\hyphenchar\font45}
\DeclareFontShape{U}{matha}{m}{n}{
  <-6> matha5 <6-7> matha6 <7-8> matha7
  <8-9> matha8 <9-10> matha9
  <10-12> matha10 <12-> matha12
  }{}

\DeclareSymbolFont{matha}{U}{matha}{m}{n}
\DeclareMathSymbol{\Lt}{3}{matha}{"CE}

\title{Generalizing K\"orner's graph entropy to graphons}

\author{Viktor Harangi}

\address{Alfr\'ed R\'enyi Institute of Mathematics, Budapest, Hungary}
\email{harangi@renyi.hu}
\thanks{During this project the first author received partial support from NRDI 
(grant KKP 138270) and from the Hungarian Academy of Sciences (J\'anos Bolyai Scholarship).}

\author{Xueyan Niu}
\author{Bo Bai}
\address{Theory Lab, Central Research Institute, 2012 Labs, Huawei Technologies Co. Ltd., Hong Kong SAR, China}
\email{\{niuxueyan3,baibo8\}@huawei.com}

\begin{document}

\maketitle

% MSC codes
%%%%%%%%%%%
% 94A29: Source coding
% 05C69: Dominating sets, independent sets, cliques
% 05C80: Random graphs

\begin{abstract}
K\"orner introduced the notion of graph entropy in 1973 as the minimal code rate of a natural coding problem where not all pairs of letters can be distinguished in the alphabet. Later it turned out that it can be expressed as the solution of a minimization problem over the so-called vertex-packing polytope.

In this paper we generalize this notion to graphons. We show that the analogous minimization problem provides an upper bound for graphon entropy. We also give a lower bound in the shape of a maximization problem. The main result of the paper is that for most graphons these two bounds actually coincide and hence precisely determine the entropy in question. Furthermore, graphon entropy has a nice connection to the fractional chromatic number and the fractional clique number. 
\end{abstract}

%%%%%%%%%%%%%%%%%%%%%%%%%%%%%%%%%%%%%%%%%%%%%%%%%%%%%%%%%%%%%%%%%%%%%%%%%%%%%%%%%%%%%%%%%
\section{Introduction}

Let $\Sigma$ be an \emph{alphabet}: a finite set, the elements of which we refer to as \emph{letters} or \emph{symbols}. Assume that a distribution $\pi$ is given on $\Sigma$: let $p_\sigma \defeq \pi(\{\sigma\})$ denote the probability of the letter $\sigma \in \Sigma$. The Shannon entropy of $\pi$ is defined as 
\[ H(\pi) = \sum_{\sigma \in \Sigma} -p_\sigma \log(p_\sigma) .\]
Shannon's classical work in source coding says that $H(\pi)$ equals the minimum \emph{code rate} for an IID sequence of letters with distribution $\pi$. That is, given an IID sequence of length $\ell$, one may encode it (in a uniquely decodable fashion) such that the expected length of the codeword is $\ell \big( H(\pi) + o(1) \big)$ as $\ell \to \infty$, and this is best possible. We can even omit the word ``expected'' in the above result if we settle for encoding only a $1-\la$ proportion (in the sense of probability) of the sequences for arbitrarily small $\la>0$.

\subsection{Distinguishable pairs and graph entropy}
Suppose now that not every pair of our alphabet's letters can be distinguished. Let $G$ be a graph with vertex set $V(G)=\Sigma$ describing which pairs are distinguishable: $x,y \in \Sigma$ can be distinguished if and only if $xy$ is an edge of $G$. Furthermore, we say that the sequences $x_1,\ldots,x_\ell$ and $y_1,\ldots,y_\ell$ are distinguishable if $x_i$ and $y_i$ are distinguishable for at least one index $i$. We wish to encode IID sequences in a way that for ``typical sequences'' it holds that distinguishable sequences are mapped to different codewords. The minimal possible code rate achievable for this problem is called \emph{graph entropy} and denoted by $H(G,\pi)$. The notion was introduced by K\"orner in \cite{korner1973}, where a relatively simple (non-asymptotic) formula was given for graph entropy. An even simpler expression was found later by Csisz\'ar, K\"orner, Lov\'asz, Marton, and Simonyi \cite{entropy_splitting}. To state their formula we need the following definition.
\begin{definition} \label{def:vp_polytope}
For a simple finite graph $G$, its \emph{vertex-packing polytope} $\VP(G)$ is defined as the convex hull of the characteristic vectors of the independent sets of $G$.

More precisely, for a subset $U$ of the vertex set $V(G)$ we write $\ind_U$ for the vector in the Euclidean space $\Rb^{V(G)}$ that has $1$ in the coordinates corresponding to $U$, and $0$ elsewhere. We refer to $\ind_U$ as the the \emph{characteristic vector} or the \emph{indicator function} of $U$. Recall that $U$ is said to be an \emph{independent set} or a \emph{stable set} if the induced subgraph $G[U]$ has no edge. By $\Ic(G)$ we denote the set of independent sets of $G$. Then 
\[ \VP(G) \defeq \conv\big( \{ \ind_J \, : \, J \in \Ic(G) \} \big) \subseteq [0,1]^{V(G)} \subset \Rb^{V(G)} .\]
\end{definition}
Given a simple finite graph and a distribution $\pi$ on the vertex set defined by the probabilities $p_x$, $x \in V(G)$, 
graph entropy can be expressed by the following simple formula:
\begin{equation} \label{eq:vp_formula_intro}
H(G,\pi) = \min_{(a_x) \in \VP(G)} \sum_{x \in V(G)} - p_x \log(a_x) .
\end{equation}
In addition to being the optimal code rate of a very natural source coding problem, graph entropy turned out to be an interesting graph theoretic tool in its own right. A surprising connection to perfect graphs was revealed by a beautiful result of Csisz\'ar, K\"orner, Lov\'asz, Marton, and Simonyi that characterizes perfect graphs in terms of graph entropy \cite{entropy_splitting}. Graph entropy has also found applications in various graph covering questions as well as in hashing and sorting problems, 
see e.g.~\cite{sorting,hashing,circuit,boolean}. 
For more background and applications, see the excellent survey papers of Simonyi \cite{survey, survey2}.

\subsection{Generalization to graphons}
Suppose now that our alphabet is not finite: say, we have a source emitting an IID sequence of uniformly random numbers from $[0,1]$. Suppose further that we do not need to distinguish between all outcomes: if $x,y \in [0,1]$ are not distinguishable, then we write $W(x,y)=0$, otherwise $W(x,y)>0$ has some positive value. As in the finite setting, we wish to encode the IID sequence in a way that distinguishable sequences are mapped to different codewords but, as before, we do not require this to hold for all sequences, only for typical sequences (say, with probability $1-\la$). The entropy of $W$ is the minimal achievable code rate for the above problem.

From this point on, we will use the language of graph limits. \emph{Graphons} are limit objects in the theory of dense graph convergence. By a graphon one normally means a symmetric measurable function $W \colon [0,1]^2 \to [0,1]$ but there is a somewhat more general treatment: given a probability space $(\Omega, \Ac, \pi)$, a graphon is a symmetric measurable function $W \colon \Omega \times \Omega \to [0,1]$.\footnote{We could restrict ourselves to the standard setting of $\Omega=[0,1]$ and $\pi$ being the Lebesgue measure but the general framework will help us to present certain examples more transparently and to explain certain phenomena better. See \cite[Chapter 13]{lovaszbook} for an overview of the basic concepts.} One can consider the analogous coding problem in this setup and define the \emph{graphon entropy} $H(W,\pi)$ as the minimal achievable code rate. (To keep the introduction concise, we need to postpone some technical details and the rigorous definitions until Section \ref{sec:prelim}.)

\begin{figure}[ht]
\centering
\includegraphics[width=3in]{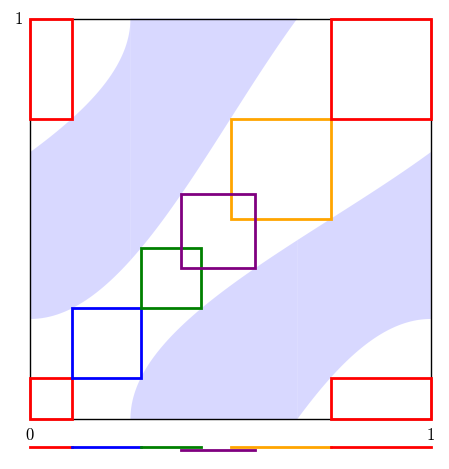} 
\caption{A graphon with independent sets: $W(x,y)$ is $1$ on the shaded part, and $0$ elsewhere; $\Omega=[0,1]$ is covered by five independent sets (depicted as colored horizontal intervals at the bottom). By definition, $J \subseteq [0,1]$ is an independent set if $W$ is almost everywhere $0$ over $J \times J$. Note that the red set $J_{\textrm{red}}$ is the union of two intervals, and hence $J_{\textrm{red}} \times J_{\textrm{red}}$ is not simply a square but the union of four rectangles.}
\label{fig:ind_sets}
\end{figure}

The notion of independent sets will play a key role in this generalized setting as well.
\begin{definition}[\cite{hladky_rocha}]
We say that $J \subseteq \Omega$ is \emph{an independent set} for $W$ if $W(x,y)=0$ for $(\pi \times \pi)$-a.e.~$(x,y)\in J \times J$. By $\Ic(W)$ we denote the set of independent sets for $W$.\footnote{When studying questions regarding independent sets of graphons, we may restrict ourselves to $0$--$1$-valued graphons 
$W \colon \Omega \times \Omega \to \{0,1\}$ often called \emph{random-free graphons}.} 

Furthermore, a graphon is said to \emph{have finite chromatic number} if there is a finite partition of $\Omega$ into independent sets: $\Omega = J_1 \cup \dots \cup J_r$ with $J_i \in \Ic(W)$ for each $i$. (See Figure \ref{fig:ind_sets} for a graphon with finite chromatic number.)
\end{definition}

Now we can define the analog of the vertex-packing polytope for graphons (which is a convex set but not always a polytope since it is defined as the convex hull of possibly infinitely many indicator functions).
\begin{definition} \label{def:vp_graphon}
By the \emph{vertex-packing set} of a graphon $W$ we mean the convex hull of the indicator functions of its independent sets: 
\begin{multline*}
\VP(W)
\defeq \conv\big( \{ \ind_J \, : \, J \in \Ic(W) \} \big) \\
= \bigg\{ \sum_{k=1}^r q_k \ind_{J_k} \, : \, r \in \Nb; q_k \geq 0; \sum_{k=1}^r q_k=1; J_k \in \Ic(W) \bigg\} ,
\end{multline*}
where the indicator function $\ind_J \colon \Omega \to \Rb$ takes value $1$ on $J$ and $0$ on $\Omega \sm J$.
\end{definition}

\subsection{Bounds for graphon entropy}
\begin{introthm} \label{thm:graphon_entropy_bounds}
Suppose that $W$ is a graphon over the probability space $(\Omega, \Ac, \pi)$. For a measurable function $a \colon \Omega \to (0,1]$ let 
\[ \fpi(a) \defeq \int_\Omega -\log(a(x)) \, \mathrm{d} \pi(x) .\]
\begin{enumerate}[(i)]
\item For any $a \in \VP(W)$ we have the upper bound\footnote{In fact, we will prove the stronger statement that $\fpi(a)$ is an upper bound for any $a = \E \ind_J$ obtained as the expectation of the indicator function of a random independent set $J\in \Ic(W)$ satisfying a technical measurability condition. See Section \ref{sec:bounds} for details.} $H(W,\pi) \leq \fpi(a)$. In other words,
\[ H(W,\pi) \leq \inf_{\VP(W)} \fpi .\]
\item We have the lower bound $H(W,\pi) \geq \fpi(b)$ whenever a measurable $b \colon \Omega \to (0,1]$ has the property that 
\[ (\ast) \quad \ds \int_J \frac{1}{b} \dd\pi \leq 1 \mbox{ for all } J \in \Ic(W) . \] 
\end{enumerate}
\end{introthm}

The main result of the paper is that we actually have matching lower and upper bounds for most graphons with finite chromatic number. 
\begin{introthm} \label{thm:graphon_matching_bounds}
Suppose that $W$ has finite chromatic number, and let $\VPcl(W)$ denote the $L_1$-closure of the vertex-packing set. Then 
\[ \inf_{\VP(W)} \fpi = \min_{\VPcl(W)} \fpi = \fpi(a) \]
for an a.e.\ unique $a \in \VPcl(W)$. If $1/a \in L_2(\pi)$ holds for this $a$, then $b=a$ satisfies condition $(\ast)$ of the lower bound, and hence 
\[ H(W,\pi)=\fpi(a) = \int_\Omega -\log(a(x)) \dd\pi(x) .\]
\end{introthm}
In particular, these theorems apply to all finite graphs (i.e., the case when $\Omega$ is finite and hence $\pi$ is an atomic measure), where the corresponding results were essentially known; see Section \ref{sec:spec_fin}. To summarize, in order to determine $H(G,\pi)$ for a finite graph $G$, it suffices to ``guess'' what the optimal $a=b \in \VP(G)$ is and simply check the condition $\sum_{x \in J} 1/b_x \leq 1$ for each $J \in \Ic(G)$ to verify optimality. 

\subsection{Connection to the fractional chromatic number}
Given a finite or infinite alphabet with distinguishable pairs, how should one choose a distribution $\pi$ on the alphabet to get the maximum possible rate for an IID sequence? This question leads to a nice connection to fractional graph parameters.

Let $G$ be a finite graph. It is well known that $\omega(G) \leq \chi(G)$, where $\omega(G)$ is the clique number (the size of the largest clique in $G$) and $\chi(G)$ is the chromatic number (the least number of colors needed for a proper coloring of $G$). It is possible to define fractional relaxations of these graph parameters (that actually coincide due to the duality of the corresponding linear programs) so that we have 
\begin{equation*}
\omega(G) \leq \omf(G) = \chif(G) \leq \chi(G) .    
\end{equation*}
The fractional chromatic number is closely related to graph entropy \cite[Lemma 4]{survey2}: 
\begin{equation} \label{eq:frac_fin}
\max_\pi H(G,\pi) = \log \chif(G) .
\end{equation}
These fractional graph parameters were defined by Hladk\'y and Rocha \cite{hladky_rocha} for graphons as well. They do not necessarily coincide any more: $\omf(W) \leq \chif(W)$. Nevertheless, we will extend \eqref{eq:frac_fin} by showing that 
\[ \log \omf(W) \leq \sup_{\pi} H(W,\pi) \leq \log \chif(W) .\]
See Section \ref{sec:frac_param} for the rigorous statements and the precise description of the probability distributions $\pi$ considered when taking the supremum.

\subsection{Set systems}
In fact, we will develop these concepts and results in the following general framework. Instead of considering the collection $\Ic(W)$ of independent sets of a graphon, we will take an arbitrary set system $\Ic$ consisting of measurable sets of the space $\Omega$. In Section \ref{sec:covering} we will define the entropy $H(\Ic,\pi)$ for any system $\Ic$ and any probability measure $\pi$ on $\Omega$, and prove our results in this context. 

We will provide plenty of examples in Section \ref{sec:examples} and show how our theorems can be applied to determine or bound the entropy. For instance, if a metric is given on $\Omega$, then it is natural to define two symbols to be distinguishable if their distance is \emph{large}. We will investigate such examples over various spaces. For a specific one, let $\Omega=S^1=\Rb \big/ \Zb$ and fix an $\alpha \in (0,1)$. We say that $x,y \in S^1$ are distinguishable if their distance (on the circle) is at least $\alpha$. Therefore, the (maximal) independent sets are \emph{arcs} of length $\alpha$: 
\[ \Ic_\alpha = \big\{ (x,x+\alpha) \, : \, x \in S^1 \big\} .\]
If we consider this system with the standard measure $\mu$ on $S^1$, then we simply have $H(\Ic_\alpha, \mu) = -\log \alpha$. Assume, instead, that we have a probability measure $\pi$ on $S^1$ that is absolutely continuous w.r.t.~$\mu$ and denote its Radon--Nikodym derivative by $g$. 
Then we have the following lower bound: 
\[ H(\Ic_\alpha, \pi) \geq 
- \int g \log g \, \mathrm{d}\mu - \log \alpha .\]
Note that the first term on the right-hand side is the \emph{differential entropy} $H(\pi)$ of $\pi$ (w.r.t.~$\mu$). We even have equality for a wide family of measures $\pi$, resulting in the neat formula $H(\Ic_\alpha, \pi) = H(\pi) - \log \alpha$. 

The above example is actually more general than it first appears. One could consider $\Omega=S^1$ with a system $\Ic$ containing arbitrary arcs. Suppose that $f \colon S^1 \to S^1$ is an \emph{orientation preserving homeomorphism} and let $\Ic=\Ic_f$ be the set of arcs of the form $\big( x, f(x) \big)$. From the classical theory of dynamical systems on the circle, we know that when the so-called \emph{rotation number} $\alpha$ is irrational, then $x \mapsto f(x)$ is \emph{semiconjugate} to the dynamical system $x \mapsto x+\alpha$. Under certain mild conditions, we can even conclude that the whole setup is equivalent to the previous example of $\Ic_\alpha$.
See Section \ref{sec:arc} for details. 

Even when the chromatic number is infinite, we may have finite entropy. For such examples, we will consider the natural setup when $\Ic$ consists of mutually independent events of a probability space $(\Omega,\pi)$. We will discuss some special cases, but it seems to be a very challenging problem to determine the entropy of these systems in full generality (Section \ref{sec:indep}).

\subsection{Graphon entropy in the literature}
Given a graphon $W$, there is a corresponding random (labeled) finite graph $G(n,W)$ on the vertex set $\{1,\ldots,n\}$. For any given $n$, $G(n,W)$ is a discrete random variable, so it is natural to consider its Shannon entropy and study its asymptotic behavior as $n \to \infty$. It turns out \cite{aldous,janson} that 
\[ \lim_{n \to \infty} \frac{ H\big( G(n,W) \big) }{\binom{n}{2}} = \iint h\big( W(x,y) \big) \, \mathrm{d}\pi(x) \mathrm{d}\pi(y) , \]
where $h(x)=-x \log x - (1-x) \log(1-x)$. In particular, this limit is $0$ for so-called \emph{random-free graphons} (taking values in $\{0,1\}$). In that case the (sub-quadratic) growth of $H\big( G(n,W) \big)$ was investigated in \cite{hatami_norine_2013}. Our entropy notion seems to have no connection to this line of research.

\subsection*{Outline of the paper}
Section \ref{sec:prelim} contains the rigorous definitions and statements. We prove the lower and upper bound in Section \ref{sec:pr_bound}, while the proof of the main result can be found in Section \ref{sec:pr_main}. The connection to the fractional chromatic number is discussed in Section \ref{sec:frac_param}. Finally, we present various examples in Section \ref{sec:examples}. 

\subsection*{Acknowledgments}
The first author thanks G\'abor Simonyi for bringing \eqref{eq:frac_fin} to his attention and suggesting that we should investigate the connection to the fractional chromatic number, as well as for further helpful comments and suggestions. The authors are also grateful to the anonymous referees: the presentation of the paper improved significantly thanks to their suggestions.

%%%%%%%%%%%%%%%%%%%%%%%%%%%%%%%%%%%%%%%%%%%%%
\section{Rigorous definitions and statements} \label{sec:prelim}

\subsection{Graphon entropy as a covering problem}

Assume that $W \colon \Omega \times \Omega \to [0,1]$ is a graphon (i.e., a measurable symmetric function) over the probability space $(\Omega,\Ac,\pi)$. We call $x,y \in \Omega$ \emph{distinguishable} if $W(x,y)>0$, while two sequences $(x_1,\ldots,x_\ell), (y_1,\ldots,y_\ell) \in \Omega^\ell$ are said to be distinguishable if there exists an index $i$ such that $x_i$ and $y_i$ are distinguishable. 

As described in the introduction, we consider the problem of encoding sequences in a way that distinguishable sequences get different codewords. This can be regarded as a covering problem due to the following simple observation (and its immediate corollary). 
\begin{proposition}
For a measurable $S \subseteq \Omega^\ell$ the following are equivalent:
\begin{enumerate}[(i)]
\item two sequences $(x_i),(y_i) \in S$ are not distinguishable almost surely (i.e., the exceptional set has $\pi^{2\ell}$-measure zero);
\item $S \subseteq (J_1 \times \ldots \times J_\ell) \cup S_0$ for some independent sets $J_1, \ldots, J_\ell \in \Ic(W)$ and a null set $S_0 \subseteq \Omega^\ell$ (i.e., $\pi^\ell(S_0)=0$).
\end{enumerate} 
\end{proposition}
\begin{proof}
By definition, a.a.~pairs of sequences in $J_1 \times \ldots \times J_\ell$ are not distinguishable, proving $(ii) \Rightarrow (i)$.

Now let $S \subseteq \Omega^\ell$ be an arbitrary measurable set. For a given $1 \leq i \leq \ell$ we define the projections $P_i \colon \Omega^\ell \to \Omega$ and $\hat{P}_i \colon \Omega^\ell \to \Omega^{\ell-1}$ as 
\[ P_i\big( (x_1, \ldots, x_\ell) \big) = x_i \mbox{ and } 
\hat{P}_i\big( (x_1, \ldots, x_\ell) \big) = (x_1,\ldots, x_{i-1},x_{i+1},\ldots,x_\ell) .\]
Let $f_i(x_i)$ be the measure of the $x_i$-section of $S$. More precisely,
\[ f_i(x_i) \defeq \pi^{\ell-1}\bigg( \hat{P}_i\big( S \cap P_i^{-1}(x_i) \big) \bigg) .\]
Furthermore, let $J_i$ denote the support of $f_i$, that is,
\[ J_i \defeq \big\{ x_i \in \Omega \, : \, f_i(x_i)>0 \big\} . \]
We claim that $J_1 \times \ldots \times J_\ell$ covers $S$ apart from a null set. Note that 
\[ S \sm (J_1 \times \ldots \times J_\ell) = \bigcup_{i=1}^\ell \big( S \cap P_i^{-1}( \Omega \sm J_i ) \big) ,\]
and for each $i$ the uncovered part has measure zero by Fubini:
\[ \pi^\ell \big( S \cap P_i^{-1}( \Omega \sm J_i ) \big) = \int_{\Omega \sm J_i} f_i  \, \mathrm{d}\pi = 0 . \]

Finally, we show that if $(i)$ holds, then each $J_i$ is an independent set, completing the proof. Assume, by contradiction, that there is a set $T \subseteq J_i^2$ consisting of distinguishable pairs $(x_i,y_i)$ such that $\pi^2(T)>0$. Then the set 
\[ \big\{ (x_1,\ldots,x_\ell,y_1,\ldots,y_\ell) \in S^2 \, : \, (x_i,y_i) \in T \big\} \]
contains pairs of distinguishable sequences and has $\pi^{2\ell}$-measure
\[ \int_T \underbrace{f_i(x_i) f_i(y_i)}_{>0}   \, \mathrm{d}\pi(x_i) \mathrm{d}\pi(y_i) > 0 \]
by Fubini, contradicting $(i)$.
\end{proof}
\begin{corollary}
For a measurable set $S \subseteq \Omega^\ell$ the following are equivalent:
\begin{enumerate}[(i)]
\item sequences in $S$ can be encoded using $N$ codewords such that sequences with the same codeword are almost surely not distinguishable;
\item $S$ can be covered, apart from a null set, by $N$ boxes of the form $J_1 \times \ldots \times J_\ell$, $J_1, \ldots, J_\ell \in \Ic(W)$.
\end{enumerate}
\end{corollary}

This means that encoding an IID sequence of length $\ell$ (with probability at least $1-\la$ and using $N$ codewords) in the desired way is equivalent to covering some set $S \subseteq \Omega^\ell$ with $\pi^\ell(S) \geq 1-\la$ by $N$ boxes (i.e., Cartesian products of independent sets). 
Note that the sole attribute of the graphon $W$ that is relevant in this problem is the set $\Ic(W)$ of independent sets. This prompts one to formulate the problem for an arbitrary collection $\Ic$ of measurable sets and investigate the analogous questions in this more general setting.

\subsection{The general covering problem} \label{sec:covering}
Let $(\Omega,\Ac,\pi)$ be a probability space and $\Ic$ an arbitrary set of measurable sets $J \in \Ac$. We will simply refer to such a (measurable) set system $\Ic$ as a \emph{system}. By an $\ell$-dimensional $\Ic$-box we mean a Cartesian product $J_1 \times \cdots \times J_\ell \subseteq \Omega^\ell$, where $J_i \in \Ic$ for each $i=1,\ldots,\ell$. Given $0<\la<1$, we wish to cover $1-\la$ proportion of $\Omega^\ell$ using as few $\Ic$-boxes as possible.

\begin{definition} \label{def:system_entropy} 
Let $N_\ell(\Ic,\pi,\la)$ be the minimum number of $\ell$-dimensional $\Ic$-boxes such that the measure of their union is at least $1-\la$ (with respect to the product measure $\pi^\ell$). We are interested in the asymptotic behaviour of the rate 
\[ \frac{1}{\ell} \log N_\ell(\Ic,\pi,\la) \mbox{ as } \ell \to \infty .\]
Since we do not necessarily know that the limit exists and is independent of $\la$, we define the \emph{lower and upper entropy} of the pair $(\Ic,\pi)$ for each $\la$ separately as follows:
\begin{align*}
\Hlow_\la(\Ic,\pi) & \defeq 
\liminf_{\ell \to \infty} \frac{1}{\ell} \log N_\ell(\Ic,\pi,\la) ; \\
\Hupp_\la(\Ic,\pi) & \defeq 
\limsup_{\ell \to \infty} \frac{1}{\ell} \log N_\ell(\Ic,\pi,\la) .
\end{align*}
We say that the \emph{entropy} $H(\Ic,\pi)$ exists and is equal to some value $0 \leq h \leq \infty$ if 
\[ \Hlow_\la(\Ic,\pi)=\Hupp_\la(\Ic,\pi)=h \mbox{ for all } 0<\la<1 .\]
For a graphon $W$, we get back the \emph{graphon entropy} $H(W, \pi)$ described in the introduction:
\[ H(W, \pi)=H(\Ic(W), \pi) .\]
\end{definition}
\begin{question}
Does the entropy $H(\Ic,\pi)$ always exist?
\end{question}

\subsection{Bounds} \label{sec:bounds}
It will be convenient to work with the following shorthand notation for the integral of $-\log(a)$.
\begin{definition} \label{def:fpi}
Given a measurable function $a \colon \Omega \to [0,\infty)$, let 
\[ \fpi(a) \defeq \int_\Omega -\log(a(x)) \dd\pi(x) = \int_\Omega \log(1/a(x)) \dd\pi(x) \]
with the usual convention that when $a(x)=0$ we have $-\log(a(x))=\log(1/a(x))=\infty$.

Note that $-\log$ is convex, and hence so is $\fpi$. Furthermore, $\fpi(a)=\infty$ if and only if $\pi\big(\{a=0\}\big) > 0$.
\end{definition}
With this notation our upper bound will read as $H(\Ic,\pi) \leq \fpi(a)$, whenever $a$ is a (finite or countable) convex combination of indicator functions of sets in $\Ic$, or more generally, whenever $a = \E \ind_J$ for some random $J \in \Ic$.

For the precise statements, we will introduce some variants of the vertex-packing polytope. 
In short, $\VP(\Ic) \defeq \conv\left( \big\{ \ind_J \, : \, J \in \Ic \big\} \right)$ consists of the (finite) convex combinations of the indicator functions $\ind_J$.\footnote{We may call $\VP(\Ic)$ the \emph{vertex-packing set} of $\Ic$.} We will also use the notation $\VPcnt(\Ic)$ for the set of countable convex combinations. Finally, for our most general bound, we will define $\VPfull(\Ic,\pi)$ as the set of all $\E \ind_J$ for which the distribution of the random set $J$ satisfies the following technical condition.
\begin{definition} \label{def:compatible}
We say that a probability measure $\beta$ over $\Ic$ is \emph{compatible} (with $\pi$) if it is defined on a $\sigma$-algebra $\Bc \subseteq \Pc(\Ic)$ in such a way that the set 
\[ \big\{ (J,x) \in \Ic \times \Omega \, : \, x \in J \big\} \]
is measurable in the product space $\Ic \times \Omega$ (w.r.t.\ the completion of the product measure $\beta \times \pi$). 
Then 
\[ F(J,x) \defeq \ind_J(x) = 
\begin{cases}
1 & \mbox{if }x \in J\\
0 & \mbox{if }x \notin J
\end{cases} \] 
defines a measurable $\Ic \times \Omega \to \{0,1\}$ function. By integrating out $J$, we get a measurable $\Omega \to [0,1]$ function:
\[ a(x) \defeq \int F(J,x) \dd\beta(J) \mbox{, or with a shorthand notation: } a \defeq \E \ind_J .\]
\end{definition}
\begin{remark}
Note that any discrete measure $\beta$ (with finite or countable support) is clearly compatible, and $\E \ind_J$ is simply a (finite or countable) convex combination in that case.
\end{remark}
\begin{definition}
We formally define the sets $\VP(\Ic)$, $\VPcnt(\Ic)$, and $\VPfull(\Ic,\pi)$ as follows:
\begin{align*}
\VP(\Ic) = & 
\big\{ \sum_{k} q_k \ind_{J_k} \, : \, J_k \in \Ic; q_k\geq 0; \sum_k q_k = 1 \mbox{ for finitely many $k$} \big\} ;\\
\VPcnt(\Ic) = & 
\big\{ \sum_{k} q_k \ind_{J_k} \, : \, J_k \in \Ic; q_k\geq 0; \sum_k q_k = 1 \mbox{ for countably many $k$} \big\} ;\\
\VPfull(\Ic,\pi) = & 
\big\{ \E_{J \sim \beta} \ind_J \, : \, \beta \mbox{ is a compatible probability distribution on $\Ic$} \} .
\end{align*}
These are convex sets satisfying 
$\VP(\Ic) \subseteq \VPcnt(\Ic) \subseteq \VPfull(\Ic,\pi)$.
Furthermore, let 
\[ K(\Ic,\pi) \defeq \big\{b \colon \Omega \to (0,1] \mbox{ measurable} 
\, : \, \int_J \frac{1}{b} \dd\pi \leq 1 \, \forall J \in \Ic \big\} .\]
\end{definition}
\begin{question}
Is there a more transparent way of defining $\VPfull$? It would also be interesting to precisely describe how $\VP$, $\VPcnt$, and $\VPfull$ are related to each other. (They often coincide, or one is the closure of another in an appropriate topology.)
\end{question}
Now we can state our lower and upper bounds. The proofs can be found in Section \ref{sec:pr_bound}.
\begin{theorem} \label{thm:all_bounds}
Let 
\begin{align} \label{eq:entropy_bounds_notations}
\begin{split}
\Hfin(\Ic,\pi) \defeq & \inf_{\VP(\Ic)} \fpi ;\\
\Hcnt(\Ic,\pi) \defeq & \inf_{\VPcnt(\Ic)} \fpi ;\\
\Hfull(\Ic,\pi) \defeq & \inf_{\VPfull(\Ic,\pi)} \fpi ;\\
\Hlow(\Ic,\pi) \defeq & \sup_{K(\Ic,\pi)} \fpi .
\end{split}
\end{align}
Then for any $0<\la<1$ we have 
\[ \Hlow(\Ic,\pi) \leq \Hlow_\la(\Ic,\pi) \leq \Hupp_\la(\Ic,\pi) \leq 
\Hfull(\Ic,\pi) \leq \Hcnt(\Ic,\pi) \leq \Hfin(\Ic,\pi) .\]
\end{theorem}

\subsection{Finite chromatic number}
A graphon $W$ is said to have finite chromatic number if $\Omega$ can be covered by finitely many independent sets of $W$. Analogously, we define the notion of chromatic number for a system $\Ic$.
\begin{definition} \label{def:chi}
Given a system $\Ic$, let $\chi(\Ic)$ be the smallest positive integer $r$ for which there exist 
\[ J_1,\ldots,J_r \in \Ic \mbox{ such that }
\pi\big( \bigcup_{k=1}^r J_k \big) = 1 ,\]
that is, $r$ sets can cover $\Omega$ apart from a null set. We write $\chi(\Ic)=\infty$ if no such finite $\Ic$-covering exists.
\end{definition}
Next we give a simple characterization of $\chi(\Ic) < \infty$. 
\begin{proposition} \label{prop:fin_chr_num}
The following are equivalent:
\begin{enumerate}[(i)]
\item $\chi(\Ic)<\infty$, that is, there exist finitely many $J_1,\ldots,J_r \in \Ic$ such that $\pi\big( \bigcup_{k=1}^r J_k \big) = 1$;
\item there exists $a \in \VP(\Ic)$ and $\delta>0$ such that $a(x) \geq \delta$ for a.e.~$x\in\Omega$;
\item there exists $a \in \VPcnt(\Ic)$ and $\delta>0$ such that $a(x) \geq \delta$ for a.e.~$x\in\Omega$. %$\, \Leftrightarrow \, \chif(\Ic)<\infty$ 
\end{enumerate}
\end{proposition}
\begin{proof}
If $\chi(\Ic)<\infty$, then using the corresponding finite covering $J_1,\ldots,J_r$ let
\[ a \defeq \frac{1}{r} \sum_{k=1}^r \ind_{J_k} \in \VP(\Ic) 
\mbox{ so that } a(x) \geq \frac{1}{r} \mbox{ for a.e.\ } x\in\Omega ,\]
proving $(i) \Rightarrow (ii)$.

The implication $(ii) \Rightarrow (iii)$ follows from the fact that $\VP(\Ic) \subseteq \VPcnt(\Ic)$.

Finally, we prove $(iii) \Rightarrow (i)$. Suppose that $a(x) \geq \delta$ for a.e.~$x\in\Omega$ for some $a = \sum_{k=1}^\infty q_k \ind_{J_k} \in \VPcnt(\Ic)$ and $\delta>0$. There exists $k_0 \in \Nb$ such that $\sum_{k>k_0} q_k < \delta/2$. It follows that 
\[ \sum_{k=1}^{k_0} q_k \ind_{J_k} \geq a - \frac{\delta}{2} 
\geq \frac{\delta}{2} \, \mbox{ a.e.} \]
In particular, $\bigcup_{k=1}^{k_0} J_k$ must have full measure and hence the chromatic number is indeed finite.
\end{proof}
\begin{remark} \label{rm:vpfull}
As Example \ref{ex:tr_copies} shows, there exists an (uncountable) $\Ic$ with infinite chromatic number for which $\VPfull(\Ic,\pi)$ contains functions with positive essential infimum.
\end{remark}

\subsection{Matching bounds}

\begin{definition}
Let $\VPcl(\Ic)$ denote the $L_1$-closure\footnote{Since $L_p$ norms (for $1\leq p < \infty$) give the same convergence notion for uniformly bounded functions over a probability space, the $L_2$-closure would be the same.}
of $\VP(\Ic)$.
\end{definition}
Now we are in a position to state our main result. See Section \ref{sec:pr_main} for the proof.
\begin{theorem} \label{thm:matching_bounds}
Suppose that $\chi(\Ic)<\infty$. Then 
\[ \inf_{\VP(\Ic)} \fpi = \min_{\VPcl(\Ic)} \fpi = \fpi(a) \]
for an a.e.\ unique $a \in \VPcl(\Ic)$. If $1/a \in L_2(\pi)$ holds for this $a$, then $H(\Ic,\pi)=\fpi(a)$. More precisely, for all $0<\la<1$ we have 
\[ \Hlow_\la(\Ic,\pi) = \Hupp_\la(\Ic,\pi) = 
\fpi(a) = \int_\Omega -\log(a(x)) \dd\pi(x) .\]
\end{theorem}
\begin{remark}
In the above theorem we may replace the condition $\chi(\Ic)<\infty$ with the weaker condition that $\VPfull(\Ic,\pi)$ contains a function $a_0$ bounded away from $0$ a.e. Then the same conclusion holds if we replace $\VP$ with $\VPfull$ in the statement. (See Remark \ref{rm:vpfull} and Example \ref{ex:tr_copies}.)

Whether the other condition (namely, $1/a \in L_2(\pi)$ for the unique minimizer $a$) is really needed, we are unsure. We seem to need the condition because Lemma \ref{lem:step3} is not true without this condition in general. However, it may be true in our specific setting $K=\VPcl(\Ic)$.
\end{remark}
\begin{question}
Do there exist measurable sets $J_k \subset \Omega$ and coefficients $q_k \geq 0$ adding up to $1$ such that 
\[ a=\sum_{k=1}^\infty q_k \ind_{J_k} \]
minimizes $\fpi$ over 
$\ds \VPcnt\big(\big\{J_k \, : \, k\in\Nb \big\}, \pi\big)$ but 
\[ \int_{J_1} \frac{1}{a(x)} \dd\pi(x) > 1 ?\]
\end{question}
%

%%%%%%%%%%%%%%%%%%%%%%%%%%%%%%%%%%%%%%%%%%%%%
\section{Proofs of the bounds} \label{sec:pr_bound}
This section contains the proof of Theorem \ref{thm:all_bounds}. Typical sets will be the key tool for proving the lower and the upper bounds.

\subsection{Typical sequences}
By a ``typical'' subset of $\Omega^\ell$ we mean a set $A$ of almost full measure (i.e., $\pi^\ell(A) \geq 1-\eps$) such that each sequence $(x_1,\ldots,x_\ell) \in A$ ``behaves averagely'' in a sense. 
\begin{lemma} \label{lem:typical}
Assume that for a measurable $a \colon \Omega \to (0,\infty)$ we have
\[ \int_\Omega |\log(a(x))| \dd \pi(x) < \infty .\]
Then for any $\eps, \de>0$ and for sufficiently large $\ell \geq L(\eps,\de)$, there exists a measurable set $A^{typ}_\ell \subset \Omega^\ell$ such that $\pi^\ell(A^{typ}_\ell) \geq 1-\eps$ and for each $(x_1,\ldots,x_\ell) \in A^{typ}_\ell$ we have
\begin{equation*}
-\fpi(a) - \delta \leq \frac{1}{\ell} \sum_{i=1}^\ell \log(a(x_i)) \leq -\fpi(a) + \delta ,
\end{equation*}
or, equivalently,
\[ \exp\big( \ell (-\fpi(a)-\delta) \big) \leq
a(x_1) \cdots a(x_\ell) \leq 
\exp\big( \ell (-\fpi(a)+\delta) \big) .\]
\end{lemma}
\begin{proof}
Given an IID sequence $x_1, x_2, \ldots$, with each $x_i$ having distribution $\pi$, and a $\pi$-integrable function $f \colon \Omega \to \Rb$, the weak law of large numbers can be applied for the IID sequence $f(x_i)$ and we get that 
\[ \frac{1}{\ell} \sum_{i=1}^\ell f(x_i) \mbox{ converges in probability to } \E_{X \sim \pi} f(X) = \int_\Omega f \dd\pi \mbox{ as } \ell \to \infty .\] 
Setting $f=\log a$, we get that the average of $\log( a(x_i) )$ converges in probability to $\ds \int_\Omega \log(a) \dd \pi = -\fpi(a)$. By definition, this means that for any given $\delta>0$ 
\[ \P\bigg( \bigg| \frac{1}{\ell} \sum_{i=1}^\ell \log(a(x_i)) + \fpi(a) \bigg| > \delta \bigg) \to 0 
\mbox{ as } \ell \to \infty ,\]
and hence this probability gets below $\eps>0$ for large enough $\ell$, and the statement of the lemma follows.  
\end{proof}

\subsection{Lower bound} \label{sec:lower}
We start with the proof of the lower bound
\[ \Hlow_\la(\Ic,\pi) \geq \Hlow(\Ic,\pi) = \sup_{K(\Ic,\pi)} \fpi .\]
That is, for any given measurable function $b \colon \Omega \to [0,1]$ with the property 
\begin{equation} \label{eq:condition}
\int_J \frac{1}{b} \, \mathrm{d} \pi \leq 1 \mbox{ for all } J \in \Ic ,
\end{equation}
we need to show that $\Hlow_\la(\Ic,\pi) \geq \fpi(b)$ for any $0<\la<1$. 
\begin{proof}
First let us assume that 
\begin{equation} \label{eq:extra_assumption}
b(x)>0 \mbox{ for all } x \in \Omega \mbox{ and } \fpi(b) < \infty .
\end{equation}
Setting $a=1/b$ we have $\fpi(b)=-\fpi(a)$. Fix $\la \in (0,1)$ and choose a positive $\eps$ less than $1-\la$ and a positive $\de$. Then, for sufficiently large $\ell$, Lemma \ref{lem:typical} provides a typical set $A^{typ}_\ell$ corresponding to $a$, $\de$, $\eps$. Setting
\[ \nu(A) \defeq \int_A a \dd\pi \]
defines a measure on $\Omega$, namely, the measure $\nu$ that is absolutely continuous w.r.t.\ $\pi$ and has Radon-Nikodym derivative $\mathrm{d}\nu/\mathrm{d}\pi=a$. Then, as a simple consequence of Fubini's theorem, we have
\begin{equation} \label{eq:RN}
\frac{\mathrm{d} \nu^\ell}{\mathrm{d} \pi^\ell}(x_1,\ldots,x_\ell) = 
\prod_{i=1}^\ell a(x_i) .
\end{equation}
Since $a(x)>0$ for all $x$, it follows that the product measures $\nu^\ell$ and $\pi^\ell$ are equivalent (i.e., absolutely continuous with respect to each other).
Also, Lemma \ref{lem:typical} says that $\eqref{eq:RN} \geq \exp\big( \ell (-\fpi(a)-\delta) \big)$ provided that $(x_1,\ldots,x_\ell) \in A^{typ}_\ell$.

For any $J \in \Ic$ we have $\nu(J) = \int_J 1/b \, \mathrm{d} \pi \leq 1$. 
It follows that an arbitrary $\Ic$-box $B=J_1 \times \cdots \times J_\ell$ has measure $\nu^\ell(B) \leq 1$. 
This means that one needs at least $\nu^\ell(A)$ $\Ic$-boxes to cover a set $A \subset \Omega^\ell$.

Now assume that the union of some $\Ic$-boxes has $\pi^\ell$-measure at least $1-\la$. Intersecting this union with $A^{typ}_\ell$, which has $\pi^\ell$-measure at least $1-\eps$, we get a set $A_\ell \subseteq A^{typ}_\ell$ with $\pi^\ell(A_\ell) \geq 1-\la-\eps>0$. However, to cover $A_\ell$ with $\Ic$-boxes, one needs at least 
\[ \nu^\ell(A_\ell) \geq 
\exp\big( \ell( \fpi(b) - \delta) \big) \pi^\ell(A_\ell) \geq 
\exp\big( \ell( \fpi(b) - \delta) \big) (1-\la-\eps) \]
boxes. Taking $\ell \to \infty$, it follows that $\Hlow_\la(\Ic,\pi)$ must be at least $\fpi(b) - \delta$. This is true for any positive $\delta$ implying that the lower entropy is at least $\fpi(b)$.

As for an arbitrary $b$, let $b_n(x) \defeq \max(b(x), 1/n)$. On the one hand, $b_n \geq b$ so \eqref{eq:condition} is satisfied by $b_n$ as well. On the other hand, $b_n(x)>0$ for all $x$, and $\fpi(b_n) \leq \fpi(1/n) = \log n < \infty$, so the extra assumption \eqref{eq:extra_assumption} holds for $b_n$. Thus, 
\begin{equation} \label{eq:bound_b_n}
\Hlow_\la(\Ic,\pi) \geq \fpi(b_n) \mbox{ for each } n.    
\end{equation}
Furthermore, $\log(1/b_n)\geq 0$ is a monotone increasing sequence converging pointwise to $\log(1/b)$. By the monotone convergence theorem it follows that the corresponding integrals converge to the integral of the limit function, so we get $\fpi(b_n) \to \fpi(b)$. This, combined with \eqref{eq:bound_b_n}, implies that $\Hlow_\la(\Ic,\pi) \geq \fpi(b)$, and the proof is complete.
\end{proof}

\subsection{Upper bound} \label{sec:upper}
The inequalities $\Hfull(\Ic,\pi) \leq \Hcnt(\Ic,\pi) \leq \Hfin(\Ic,\pi)$ follow readily from the fact that $\VPfull(\Ic,\pi) \supseteq \VPcnt(\Ic) \supseteq \VP(\Ic)$.
Therefore, it remains to show that $\Hupp_\la(\Ic,\pi) \leq \Hfull(\Ic,\pi)$, that is, for a random set $J \in \Ic$ with a compatible distribution $\beta$ (recall Definition \ref{def:compatible}) it holds that 
\[ \Hupp_\la(\Ic,\pi) \leq \fpi(a) \mbox{, where } a= \E \ind_J . \]
\begin{proof}
The idea is to use random $\Ic$-boxes to cover a large portion of $\Omega^\ell$. In essence, the probability that a typical sequence is covered by a random $\Ic$-box $J_1 \times \cdots \times J_\ell$ (with the $J_i$'s drawn independently from the given distribution $\beta$) is $\exp\big( \ell(-\fpi(a)+o(1)) \big)$. It will follow that a large portion of the typical set of sequences (and hence a large portion of the whole space $\Omega^\ell$) can be covered using $\exp\big( \ell(\fpi(a)+o(1)) \big)$ boxes. Therefore the entropy in question is indeed at most $\fpi(a)$. Next we give the formal argument.

For a fixed $x \in \Omega$ and a $\beta$-random $J$, the event $\{x \in J\}$ is measurable and its probability (i.e., its $\beta$-measure) is equal to $a(x)$ for a.e.~$x$. 

Let $J_1, \ldots, J_\ell$ be independent, each with distribution $\beta$. Then for the random $\Ic$-box $J_1 \times \cdots \times J_\ell$ we have
\begin{equation} \label{eq:pt_box}
\P(x \in J_1 \times \cdots \times J_\ell) = \prod_{i=1}^\ell \P(x_i \in J_i) 
= a(x_1) \cdots a(x_\ell) 
\end{equation}
for a.e.\ sequence $x=(x_1,\ldots,x_\ell) \in \Omega^\ell$. 

Now for $\de>0$ and $0<\eps<\la$ let $A^{typ}_\ell$ denote the typical set provided by Lemma \ref{lem:typical}. By removing a zero-measure set from $A^{typ}_\ell$ we may assume that \eqref{eq:pt_box} holds for all $x \in A^{typ}_\ell$. Then the probability that any fixed typical sequence $x \in A^{typ}_\ell$ lies in the random $\Ic$-box $J_1 \times \cdots \times J_\ell$ is at least 
\[ \exp\big( \ell(-\fpi(a)-\de) \big) .\]

Let us take $M$ independent copies of such random $\Ic$-boxes for $M=\exp\big( \ell(\fpi(a)+2\de) \big)$, and let $S$ denote their union. Then any fixed $x \in A^{typ}_\ell$ is \emph{not} covered by $S$ with probability at most
\[ \bigg( 1 - \exp\big( \ell(-\fpi(a)-\de) \big) \bigg)^M \leq 
\left( \frac{1}{e} \right)^{\exp(\ell \de)} .\]
Consequently, for any given $x \in A^{typ}_\ell$ we have 
\[ \P( x \in S ) = 1-o(1) \mbox{ as } \ell\to\infty .\]
It follows that
\[ \E \, \pi^\ell(S) \geq (1-o(1)) \pi^\ell(A^{typ}_\ell) \geq (1-o(1))(1-\eps) ,\]
which is larger than $1-\la$ for large enough $\ell$. (Here the expectation is with respect to the randomness of the $M$ boxes.) We conclude that there exist $M$ $\Ic$-boxes whose union $S$ has measure $\pi^\ell(S) \geq 1-\la$ provided that $\ell$ is sufficiently large. By Definition \ref{def:system_entropy} this means that $N_\ell(\Ic,\pi,\la) \leq M = \exp\big( \ell(\fpi(a)+2\de) \big)$ implying $\Hupp_\la(\Ic,\pi) \leq \fpi(a) + 2\de$. This holds for all $\de>0$, and hence $\Hupp_\la(\Ic,\pi) \leq \fpi(a)$, completing the proof of the upper bound (and Theorem \ref{thm:all_bounds}).
\end{proof}

\section{Proof of the main result} \label{sec:pr_main}
In this section we prove Theorem \ref{thm:matching_bounds} which states that our (weakest) upper bound $\Hfin(\Ic,\pi)=\inf_{\VP(\Ic)} \fpi$ actually coincides with our lower bound under fairly mild conditions, determining the precise value of $H(\Ic,\pi)$ in those cases. Our proof strategy will be as follows.
\begin{itemize}
\item In the formula $\ds \inf_{\VP(\Ic)} \fpi$ one may replace $\VP(\Ic)$ with its closure $\VPcl(\Ic)$. 
\item Over $\VPcl(\Ic)$ the infimum is attained so it is actually a minimum: $\ds \min_{\VPcl(\Ic)} \fpi$. 
\item The unique minimizer $a \in \VPcl(\Ic)$ satisfies condition \eqref{eq:condition} so $\fpi(a)$ is a lower bound as well, implying $H(\Ic,\pi)=\fpi(a)$. 
\end{itemize}
First we prove a series of standalone lemmas that will be needed when we put all the ingredients together and present the proof of Theorem \ref{thm:matching_bounds} at the end of the section.

\subsection{Taking the closure}
\begin{lemma} \label{lem:step1}
Suppose that $\chi(\Ic) < \infty$ for a system $\Ic$. Then 
\[ \inf_{\VPcl(\Ic)} \fpi = \inf_{\VP(\Ic)} \fpi .\]
\end{lemma}
\begin{proof}
The inequality $\leq$ is trivial. For the other inequality, let $a_0 \in \VPcl(\Ic)$ be arbitrary. Since $\chi(\Ic) < \infty$, there exists $a_1 \in \VP(\Ic)$ and $\delta>0$ such that $a_1(x) \geq \delta$ for a.e.~$x\in\Omega$. For $0 \leq t \leq 1$ we set 
\[ a_t \defeq (1-t)a_0 + t a_1 .\]
Due to the convexity of $\fpi$, it follows that 
\[ \fpi(a_t) \leq (1-t)\fpi(a_0) + t \fpi(a_1) = 
\fpi(a_0) + t \big( \fpi(a_1) - \fpi(a_0) \big) .\]
Since $\fpi(a_1) \leq -\log(\delta) < \infty$ and $\fpi(a_0) \geq 0$, it follows that for any $\eps>0$ we can choose a fixed $t>0$ such that 
\[ \fpi(a_t) \leq \fpi(a_0) + \eps .\]
Since $a_0$ lies in the closure $\VPcl(\Ic)$, there exist $b_1, b_2, \ldots \in \VP(\Ic)$ such that $b_n \to a_0$ in $L_1$. Now we take the following convex combinations (with the same fixed $t$):
\[ b_{n,t} \defeq (1-t) b_n + t a_1 \in \VP(\Ic) .\]
As $n\to \infty$, $b_{n,t} \to a_t$ in $L_1$ and all these functions are uniformly bounded away from $0$. (Indeed, their essential infimum is $\geq t \delta$.) For such functions $\fpi$ is clearly Lipschitz w.r.t.\ the $L_1$ distance. It follows that 
\[ \fpi(b_{n,t}) \to \fpi(a_t) \leq \fpi(a_0) + \eps ,\]
and hence 
\[ \inf_{\VP(\Ic)} \fpi \leq \fpi(a_0) + \eps .\]
Since this is true for any $a_0 \in \VPcl(\Ic)$ and any $\eps>0$, the proof is complete. 
\end{proof}

\subsection{The minimum is attained}
Since $\VPcl(\Ic)$ is not necessarily compact, we need to do some work to conclude that the minimum is attained over $\VPcl(\Ic)$. We claim that if we have a sequence $a_n \in \VP(\Ic)$ such that $\fpi(a_n)$ converges to the infimum, then $a_n$ is a Cauchy sequence and hence convergent in the $L_2$ metric. The next lemma proves this in a more general setting. Note that it applies to our situation as $-\log$ is ``homogeneously strictly convex'' on $(0,1]$ in the sense that its second derivative ($1/t^2$) is bounded away from $0$ on $(0,1]$. 
\begin{lemma} \label{lem:step2}
Let $T$ be a (finite or infinite) interval of $\Rb$ and $f \colon T \to \Rb$ be a (convex) function such that $f''(t) \geq \delta > 0$ for all $t \in T$. Let $K$ be a convex set of measurable $\Omega \to T$ functions and for $a \in K$ let 
\[ \kappa_f(a) \defeq \int_\Omega f(a(x)) \dd\pi(x) .\]
Suppose that the infimum $s \defeq \inf_{a\in K} \kappa_f(a)$ is finite, and $\kappa_f(a_1)$ and $\kappa_f(a_2)$ are $\eps$-close to the infimum for some $\eps>0$, that is, $\kappa_f(a_1),\kappa_f(a_2) \leq s+\eps$. Then the $L_2$ distance (and hence the $L_1$ distance) of $a_1$ and $a_2$ can be bounded as follows: 
\[  \| a_1 - a_2 \|_1 \leq \| a_1 - a_2 \|_2 \leq 
\sqrt{\frac{8\eps}{\delta}} .\]
In particular, if we have a sequence $a_n \in K$ such that $\kappa_f a_n \to s$ as $n \to \infty$, then $a_n$ must be a Cauchy sequence w.r.t.\ $L_2$ norm. Since $L_2$ spaces are complete, it follows that there is a unique limiting function $\hat{a}$ for such sequences and it is meaningful to write 
\[ \lim_{\kappa_f a \to s} a = \hat{a} .\]
\end{lemma}
\begin{proof}
Since $K$ is convex, $(a_1+a_2)/2 \in K$ and hence 
\[ \kappa_f\left( \frac{a_1+a_2}{2} \right) \geq s .\]
Then
\begin{multline*}
\eps \geq \frac{\kappa_f a_1 + \kappa_f a_2}{2} - \kappa_f\left( \frac{a_1+a_2}{2} \right) = 
\int_\Omega \frac{f(a_1(x)) + f(a_2(x))}{2}- f\left( \frac{a_1(x)+a_2(x)}{2} \right) \dd\pi(x) \geq \\
\int_\Omega \frac{\delta}{8} \big( a_1(x)-a_2(x) \big)^2 \dd\pi(x) = \frac{\delta}{8} \| a_1-a_2 \|_2^2 .
\end{multline*}
Since $\pi$ is a probability measure, it follows that 
\[ \| a_1 - a_2 \|_1 \leq \| a_1 - a_2 \|_2 \leq \sqrt{\frac{8\eps}{\delta}} .\]
\end{proof}

\subsection{The minimizer}
As a final ingredient, we show that the minimizer of $\fpi$ over $K$ must lie in the so-called antiblocker of $K$ under mild assumptions.
\begin{lemma} \label{lem:step3}
Suppose that $K$ is a convex set of measurable $\Omega \to [0,1]$ functions and $a \in K$ with $1/a \in L_2(\pi)$, i.e., 
\[ \int_\Omega \frac{1}{(a(x))^2} \dd \pi(x) < \infty .\]
If $a \in K$ minimizes $\fpi$ over $K$, i.e., 
\[ \fpi(a) = \inf_K \fpi ,\]
then $1/a$ lies in the antiblocker of $K$, i.e., 
\[ \int_\Omega \frac{b(x)}{a(x)} \dd\pi(x) \leq 1 \mbox{ for all } b \in K .\]
\end{lemma}
\begin{proof}
Under stronger conditions we could argue that the directional derivative of $\fpi$ at $a$ in the direction $b-a$ is 
\[ \int_\Omega \frac{-1}{a} (b-a) \dd\pi ,\]
which is nonnegative when $\fpi$ has a minimum at $a$, and the statement of the lemma would follow. We need to be more careful under our mild conditions.

Let $b \in K$ and consider the point $a_t \defeq (1-t)a+tb=a + t (b-a)$ as $t \to 0+$, that is, a point on the line connecting $a$ and $b$, and approaching $a$. Since $a_t \in K$, we have $\fpi(a_t) \geq \fpi(a)$ for any $0<t<1$. 

We will use the bound $-\log(1+x) \leq -x+x^2$ that holds for any $x \in [-1/2,\infty)$. For $t \leq 1/2$ we get that 
\begin{multline*}
0 \leq \fpi(a_t) - \fpi(a) = \int_\Omega -\log\left( 1 + t \big( \frac{b(x)}{a(x)} - 1 \big) \right) \dd\pi(x) \leq \\
-t \bigg( \int_\Omega \frac{b(x)}{a(x)} \dd\pi(x) - 1 \bigg) + t^2 \int_\Omega \bigg( \frac{b(x)}{a(x)} - 1 \bigg)^2 \dd\pi(x) .
\end{multline*}
The second term of the right-hand side is at most $Ct^2$ for some $C<\infty$ (because $b \in L_\infty(\pi)$ and $1/a \in L_2(\pi)$). Taking $t \to 0$, it follows that 
\[ \int_\Omega \frac{b(x)}{a(x)} \dd\pi(x) \leq 1 \mbox{ for all } b \in K ,\]
as claimed.
\end{proof}

\subsection{Proof of Theorem \ref{thm:matching_bounds}}
By Lemma \ref{lem:step1} we have 
\[
s \defeq \inf_{\VP(\Ic)} \fpi = \inf_{\VPcl(\Ic)} \fpi .
\]
Using Lemma \ref{lem:step2} with $K=\VP(\Ic)$ and $f=-\log$ tells us that if we have a sequence $a_n \in \VP(\Ic)$ such that $\fpi(a_n)$ converges to the infimum $s$, then $a_n$ is a Cauchy sequence and hence convergent in the $L_2$ metric. So $a_n$ converges in $L_2$ to a limit point $a \in \VPcl(\Ic)$. We would like to conclude that $\fpi(a)=s$ and hence the minimum is attained at $a$. We need to be careful as $\fpi$ is not quite continuous w.r.t.\ the $L_2$ metric. Instead, notice that $L_2$ convergence implies that a subsequence converges almost surely to $a$. So we may assume $a_n(x) \to a(x)$ for a.e.~$x$. For measurable functions $a \colon \Omega \to [0,1]$, $x \mapsto -\log(a(x))$ is a measurable $\Omega \to [0,+\infty]$ function. Therefore Fatou's lemma can be applied for the sequence $-\log(a_n(x))$ and we get 
\[ \fpi(a) \leq \liminf_{n \to \infty} \fpi(a_n) = s .\]
Since $a \in \VPcl(\Ic)$, $\fpi(a) \geq s$ holds as well. Thus $\fpi(a)=s$ as claimed. 

It also follows that the minimizer must be a.e.\ unique. Indeed, simply apply Lemma \ref{lem:step2} with $\eps=0$. 

Finally, Lemma \ref{lem:step3} shows that if $1/a \in L_2(\pi)$ for this unique minimizer $a$, then $a$ must lie in the antiblocker. In particular, for any $J \in \Ic$, the inner product of $1/a$ and $\ind_J \in K = \VP(\Ic)$ is at most $1$, which confirms that $a$ satisfies condition \eqref{eq:condition}, and hence $\fpi(a)$ is a lower bound (as well as an upper bound) for the entropy, and the proof is complete.

\subsection{Finite graphs} \label{sec:spec_fin}
If we take a $0$--$1$-valued graphon $W$ over some finite $\Omega$, then we get back the case of finite graphs, and our general theorem turns into the following result. This has been known, although not stated in this exact form. (It follows easily from \cite[Corollary 6]{entropy_splitting} using the fact that $\VP(G)$ and the so-called fractional vertex-packing polytope of the complement graph form an antiblocking pair.)
\begin{corollary} \label{cor:finite}
Let $G$ be a finite simple graph. Suppose that $p_x>0$ for all $x \in V(G)$ with $\sum_x p_x=1$, defining a probability distribution $\pi$ on $V(G)$. Recall that $\VP(G)$ stands for the vertex-packing polytope; see Definition \ref{def:vp_polytope}. Furthermore, let $K(G,\pi)$ denote the set of points $a=(a_x) \in \Rb^{V(G)}$ satisfying the condition
\begin{equation} \label{eq:fin_condition}
\sum_{x \in J} \frac{p_x}{a_x} \leq 1 \mbox{ for each independent set } J \in \Ic(G) .
\end{equation}
Then both $\VP(G)$ and $K(G,\pi)$ are convex sets, and their intersection is a single point which is both the minimizer of the function $\fpi \colon a \mapsto \sum - p_x \log(a_x)$ over $\VP(G)$ and the maximizer of the same function over $K(G,\pi)$. This minimum/maximum value is equal to the graph entropy $H(G,\pi)$. In other words, to determine $H(G,\pi)$, it suffices to find a point $a \in \VP(G)$ that satisfies condition \eqref{eq:fin_condition}.
\end{corollary}

We include a simple example demonstrating how to use the optimality condition \eqref{eq:fin_condition}. 
\begin{proposition}
Let $G = C_{2n+1}$ be the cycle $x_0 x_1 \cdots x_{2n}$ of length $2n+1$. Suppose that 
\[ \pi(\{x_k,x_{k+1}\}) = p_{x_k} + p_{x_{k+1}} \leq 1/n  \mbox{ for each } k=0,\ldots,2n ,\]
where $x_{2n+1}=x_0$. Then we have the following simple formula for graph entropy:
\[ H(C_{2n+1}, \pi) = H(\pi) - \log(n) .\]
\end{proposition}
\begin{proof}
We will use indices modulo $2n+1$. The maximum-size independent sets are the following:
\[ J_k = \big\{ x_{k+2}, x_{k+4}, \ldots, x_{k+2n} \big\}
;\ k=0,1, \ldots, 2n .\]
For $n \geq 4$ there are other maximal independent sets but, under our assumptions on $\pi$, only these maximum-size independent sets will be needed. 

Let 
\[ q_k \defeq 1 - n \big( p_{x_k} + p_{x_{k+1}} \big) \geq 0 
\mbox{ for each } k=0,1,\ldots,2n .\] 
Note that the sum of $q_k$ is $1$, therefore 
\[ a \defeq \sum_{k=0}^{2n} q_k \ind_{J_k} \in \VP(C_{2n+1}).\]
We get that 
\[ a_{x_k} = q_{k-2} + q_{k-4} + \cdots + q_{k-2n} = n - n 
\underbrace{\sum_{i=1}^{2n} p_{x_{k-i}}}_{=1-p_{x_k}} = n p_{x_k} .\]
Thus, $p_{x_k}/a_{x_k}=1/n$ for each $k$, and hence \eqref{eq:fin_condition} holds for any set $J$ of size at most $n$, in particular, it holds for all independent sets of $C_{2n+1}$. Consequently, due to Corollary \ref{cor:finite}, $a$ is the ``optimal'' point in $\VP(C_{2n+1})$, and hence 
\[ H(C_{2n+1}, \pi) = \fpi(a) = \sum_{k=0}^{2n} -p_{x_k} \log\big(n p_{x_k}\big) = H(\pi) - \log(n) . \]
\end{proof}
%

%%%%%%%%%%%%%%%%%%%%%%%%%%%%%%%%%%%%%%%%%%%%%%%%%%%%%%%%%%%%%%%%%%%%%%%%%%%%%%%%%%%%%%%%%%%%%%%%%%%%%%%%%%
\section{The fractional clique and chromatic number} \label{sec:frac_param}

Let $(\Omega,\Ac, \mu)$ be a probability space and $\Ic$ a system. We use $\mu$ here (instead of $\pi$) because in this section we will consider the entropy $H(\Ic,\pi)$ for various probability measures $\pi$ (that are absolutely continuous w.r.t.\ $\mu$). The only role $\mu$ plays is to determine the zero-measure sets; it could be replaced with any equivalent probability measure. 
When defining the parameters below, only the set of zero-measure sets is relevant. For that reason, we will be somewhat sloppy and omit $\mu$ from the notations.\footnote{We did the same when defining the chromatic number $\chi(\Ic)$ in Definition \ref{def:chi}, where $\Omega$ was required to be covered with $\Ic$-sets only up to a set of measure zero.} 
\begin{definition}
We define the \emph{fractional chromatic number} of the system $\Ic$ (and the measure $\mu$) as 
\[ \chif(\Ic) \defeq \inf \left\{ \sum_{J \in \Ic} c(J) \, : \, 
c \colon \Ic \to [0,1] \mbox{ s.t.~for a.e.~} x \in \Omega: \, 
\sum_{J \, : \, x \in J \in \Ic} c(J) \geq 1 \right\} .\]
Note that the only place where the measure is used is hidden in the ``a.e.''\ statement which only hinges on the zero-measure sets. It is easy to see that 
\[ \frac{1}{\chif(\Ic)} = \sup \left\{ t \, : \, 
\exists a \in \VPcnt(\Ic) \mbox{ s.t.~} a \geq t \mbox{ a.e.} \right\} .\]
\end{definition}
\begin{definition}
We define the \emph{fractional clique number} of the system $\Ic$ (and the measure $\mu$) as 
\[ \omf(\Ic) \defeq \sup \left\{ \int_\Omega b \, : \, 
b \colon \Omega \to [0,1] \mbox{ measurable s.t.~for each } J \in \Ic: \, 
\int_J b \leq 1 \right\} .\]
At a first glance, it seems that this definition does depend on the measure (through the integrals). However, one could use the following equivalent formulation:
\[ \omf(\Ic) = \sup \big\{ \nu(\Omega) \, : \, 
\nu \ll \mu \mbox{ is a probability measure s.t. for each } J \in \Ic: \, 
\nu(J) \leq 1 \big\} ,\]
which, in fact, shows that the parameter only depends on the class of probability measures that are absolutely continuous w.r.t.\ $\mu$, which, in turn, is determined by the zero-measure sets. 
\end{definition}
For $\Ic=\Ic(W)$ we get back the graphon analogues of these parameters, see \cite{hladky_rocha}.

\begin{proposition}
For any probability measure $\pi \ll \mu$ we have 
\[ \Hcnt(\Ic,\pi) \leq \log \chif(\Ic). \]
\end{proposition}
\begin{proof}
By definition, for any $t>\chif(\Ic)$ there exists $a \in \VPcnt(\Ic)$ such that $a \geq 1/t$ a.e. Then for any $\pi \ll \mu$ we get 
\[ \Hcnt(\Ic,\pi) \leq \fpi(a) \leq \log(t) ,\] 
and the inequality follows.
\end{proof}
\begin{proposition}
We have 
\[ \sup_{\pi \ll \mu} \Hlow(\Ic,\pi) \geq \log \omf(\Ic). \]
\end{proposition}
\begin{proof}
By definition, for any $t<\omf(\Ic)$ there exists $\nu \ll \mu$ with $\nu(\Omega)=t$ and $\nu(J) \leq 1$ for each $J \in \Ic$. It follows that for the probability measure $\pi \defeq \frac{1}{t}\nu$ the constant function $\frac{1}{t} \ind_\Omega$ lies in $K(\Ic,\pi)$, and hence $\Hlow(\Ic,\pi) \geq \fpi\big( \frac{1}{t} \big) = \log t$.
\end{proof}
Combining these inequalities yields the following result when $\omf(\Ic)$ and $\chif(\Ic)$ coincide.
\begin{theorem}
Suppose that $\omf(\Ic)=\chif(\Ic)$ for a system $\Ic$ (and a measure $\mu)$. Then 
\[ \sup_{\pi} H(\Ic,\pi) = \log \chif(\Ic) ,\]
where the supremum is taken over all probability measures $\pi \ll \mu$.
\end{theorem}
Note that for finite graphs $G$ we always have $\omf(G)=\chif(G)$, and hence we get back \cite[Lemma 4]{survey2}, see \eqref{eq:frac_fin}.
It would be interesting to study this supremum for systems $\Ic$ with $\omf(\Ic)<\chif(\Ic)$. For instance, Hladký and Rocha noted in \cite{hladky_rocha} that a construction of Leader \cite{leader} can be turned into a graphon with different fractional clique and chromatic numbers.

%%%%%%%%%%%%%%%%%%%%%%%%%%%%%%%%%%%%%%%%%%%%%%%%%%%%%%%%%%%%%%
\section{Examples} \label{sec:examples}
This section is devoted to the analysis of various families of examples.
\subsection{Distance-based distinguishability}

It is very natural to consider examples where symbols can be distinguished if their distance, in some sense, is sufficiently large. 

\begin{example} \label{ex:dist_based}
Suppose that we have a measurable symmetric function $d \colon \Omega \times \Omega \to [0,\infty)$. Intuitively, $d$ expresses some kind of distance between pairs of symbols, and we say that two symbols are distinguishable if their distance is at least some given constant $c>0$. In other words, we consider the following graphon:
\begin{equation} \label{eq:distance_graphon}
W(x,y)= 
\begin{cases}
1 & \mbox{if } d(x,y)>c \\
0 & \mbox{otherwise.}
\end{cases} 
\end{equation}
In many natural examples $d$ is actually a metric on $\Omega$, in which case $\Ic(W)$ consists of sets of diameter at most $c$.
\end{example}
\begin{proposition}[Example \ref{ex:dist_based} for ``homogeneous'' metric measure spaces]
Suppose that the metric measure space $(\Omega,\pi,d)$ has the following properties for some fixed $c>0$.
\begin{itemize}
\item It is homogeneous in the sense that closed balls $\overline{B}(x,c/2)$ of fixed radius $c/2$ have the same measure:
\[ \pi\left( \overline{B}(x,c/2) \right) =  m_c 
\quad (\forall x \in \Omega) .\]
\item The ball has the largest measure among sets of diameter $c$: 
\[ \mbox{if } \diam(A) \leq c \mbox{, then } \pi(A) \leq m_c . \] 
\end{itemize}
Then the entropy of the graphon \eqref{eq:distance_graphon} is 
\[ H(W,\pi) = -\log m_c .\]
\end{proposition}
\begin{proof}
The constant function $a(x)=m_c$ has $\fpi(a)=-\log m_c$. Since each independent set $J$ has essential diameter at most $c$, its measure $\pi(J)$ is at most $m_c$. Therefore 
\[ \int_J \frac{1}{a} \, \mathrm{d} \pi \leq \frac{\pi(J)}{m_c} \leq 1 ,\] 
and hence $\fpi(a)$ is a lower bound for the entropy. To see that it is also an upper bound, we need to show that $a \in \VPfull(\Ic(W))$. Take a $\pi$-random point $x \in \Omega$ and the corresponding random ball $J=\overline{B}(x,c/2)$. For this random independent set $J$ we clearly have $\E \ind_J=a$.
\end{proof}

A simple special case is the circle $S^1$.
\begin{corollary} \label{cor:S1}
\label{cor:S1_basic} 
Consider Example \ref{ex:dist_based} with the following setup.
\begin{itemize}
\item Let $\ds \Omega=S^1=\Rb \big/ \Zb = [0,1] \big/_{0 \sim 1}$ with the standard measure $\mu$.
\item Let $d(x,y)$ be the length of the shorter arc between $x$ and $y$: 
\[ d(x,y)=\min(|x-y|,1-|x-y|) .\]
\end{itemize}
That is, for a fixed $0<c<1/2$ we consider the following graphon $W \colon [0,1]\times [0,1] \to \{0,1\}$:
\[ W(x,y) =
\begin{cases}
1 & \mbox{if } c< |x-y| < 1-c ;\\
0 & \mbox{otherwise.} 
\end{cases}
\]
Then $H(W,\mu)=-\log c = \log(1/c)$.
\end{corollary}
\begin{figure}[ht]
\centering
\includegraphics[width=2in]{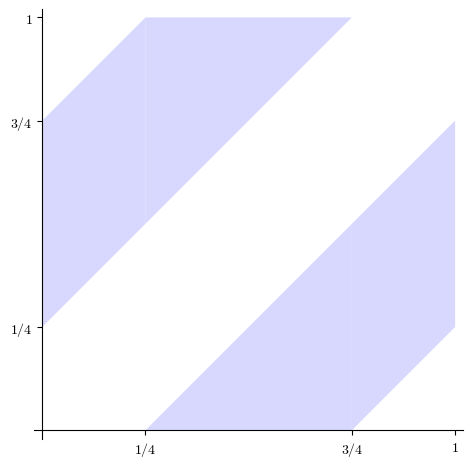} 
\hspace{0.5in} 
\includegraphics[width=2in]{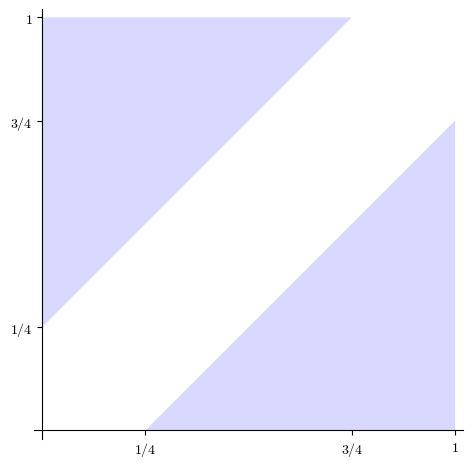} 
\caption{The graphons of Corollary \ref{cor:S1} (left) and Proposition \ref{prop:interval} (right) for $c=1/4$; $W(x,y)$ is $1$ on the shaded part, and $0$ elsewhere}
\label{fig:S1_vs_interval}
\end{figure}
The next example does not have the property that the ball around each point has the same measure.
\begin{proposition} \label{prop:interval}
Consider Example \ref{ex:dist_based} with the following setup.
\begin{itemize}
\item Let $\Omega=[0,1]$ with the Lebesgue measure $\mu$.
\item Let $d(x,y)=|x-y|$. 
\end{itemize}
That is, for a fixed $0 < c \leq 1$ we consider the following graphon: 
\[ W(x,y) =
\begin{cases}
1 & \mbox{if } |x-y|>c;\\
0 & \mbox{otherwise.} 
\end{cases}
\]
If $n$ is the unique positive integer such that $\frac{1}{n+1} < c \leq \frac{1}{n}$. Then 
\[ H(W,\mu) = n(n+1) \bigg( \left(\frac{1}{n}-c\right)\log(n+1) + 
\left(c-\frac{1}{n+1} \right) \log n \bigg) .\]
\end{proposition}
See Figure \ref{fig:S1_vs_interval} to compare the graphons in Corollary \ref{cor:S1} and Proposition \ref{prop:interval}. 
\begin{proof}
We will use the following coefficients:
\[ q_k \defeq \frac{n-k+1}{n(n+1)}\ ; \quad k=1,2,\ldots,n .\]
Note that $q_1+\cdots+q_n=1/2$. Therefore 
\[ a \defeq \sum_{k=1}^n q_k \bigg( 
\ind_{[(k-1)c,kc]} + \ind_{[1-kc,1-(k-1)c]}  \bigg) \in \VP( \Ic(W) ) .\]
We set $r=1-nc$. Note that $0 \leq r < c$. Easy calculations show that 
\[ a(x) = 
\begin{cases}
1/(n+1) & \mbox{if } 
x \in [0,r) \cup (c,c+r) \cup (2c,2c+r) \cup \cdots \cup (nc,1];\\
1/n & \mbox{otherwise.} 
\end{cases}
\]
We claim that $a$ satisfies condition \eqref{eq:condition} of the lower bound, which would imply $H(W,\mu)=\fpi(a)$, completing the proof. It suffices to check the condition for the independent set $J=[t,t+c]$ for any given $t \in [0,1-c]$:
\[ \int_t^{t+c} 1/a(x) \, \mathrm{d}x = r(n+1) + (c-r)n = r+cn = 1 .\]
\end{proof}

\subsection{Arc systems} \label{sec:arc}
Now we consider systems consisting of arcs of the circle $S^1$. For $x,y \in S^1$ we write $(x,y)$ for the open arc that starts at $x$ and goes in the positive direction ending at $y$.
\begin{example}
Let $\Omega = S^1$ with the standard measure $\mu$. For an \emph{orientation-preserving homeomorphism} $f \colon S^1 \to S^1$ let 
\[ \Ic_f \defeq \left\{ \big( x, f(x) \big) \, : \, x \in S^1 \right\} \]
Note that there is a corresponding graphon $W_f \colon S^1 \times S^1 \to \{0,1\}$ for which $W_f(x,y)=0$ if and only if $y$ lies on the arc $\big(x, f(x) \big)$. It is easy to see that $\Ic_f$ consists of the maximal independent sets (modulo null sets) of $W_f$. 

In fact, Figure \ref{fig:ind_sets} in the introduction depicts such a graphon $W_f$. Specifically, after identifying $S^1$ with $[0,1)$, we used the function 
\[ f(x)=\bigg\{ \frac{3}{4}-\frac{1}{2}\cos(\pi x) \bigg\} . \]
\end{example}

Considered as a dynamical system on the circle, the key invariant of $f$ is the so-called \emph{rotation number} $\alpha \in (0,1)$. Intuitively, $\alpha$ is the average rotation along an orbit of $f$.
\begin{proposition}
Suppose that the rotation number $\al$ of $f$ is irrational. Then there exists $G \colon S^1 \to S^1$ such that 
\begin{equation} \label{eq:conjugate}
G(f(x)) = G(x) + \alpha \mbox{ for all } x \in S^1 \quad \mbox{(addition is meant modulo $1$)},
\end{equation}
and we have the following lower bound for the entropy of the system $\Ic_f$:
\[ H(\Ic_f,\mu) \geq \int \log(G') \,\mathrm{d}\mu - \log \alpha .\]
\end{proposition}
\begin{proof}
A classical result of Poincar\'e says that if the rotation number $\alpha$ is irrational, then $f$ is semiconjugate to $x \mapsto x+\alpha$: there exists a ``monotone'' continuous map $G \colon S^1 \to S^1$ satisfying \eqref{eq:conjugate}. Such a monotone function is a.e.\ differentiable, its derivative $G' \geq 0$ is measurable and satisfies 
\[ \int_{x_1}^{x_2} G' \,\mathrm{d}\mu \leq G(x_2) - G(x_1) .\]
It follows that the function $b \defeq \alpha/G'$ satisfies condition \eqref{eq:condition} of the lower bound:
\[ \int_{x}^{f(x)} \frac{1}{b} \,\mathrm{d}\mu = 
\frac{1}{\alpha} \int_{x}^{f(x)} G' \,\mathrm{d}\mu \leq 
\frac{1}{\alpha} \bigg( \underbrace{G(f(x))}_{=G(x)+\alpha} - G(x) \bigg) 
= 1 .\]
Therefore, 
\[  H(\Ic_f,\mu) \geq \varphi_\mu(b) = \int \log(G') \,\mathrm{d}\mu - \log \alpha .\]
\end{proof}
When $G \colon S^1 \to S^1$ is a bijection (which, according to Denjoy's theorem, can be assumed when $f$ is, for instance, a $C^2$ map), then the pair $(\Ic_f,\mu)$ can be considered to be equivalent to a pair $(\Ic_\alpha,\pi)$ for some measure $\pi$, where 
\[ \Ic_\alpha \defeq \left\{ \big( x, x+\alpha \big) \, : \, x \in S^1 \right\} .\]
The next example describes this equivalent situation for the case when $\pi$ is absolutely continuous w.r.t.\ $\mu$.
\begin{proposition} \label{prop:arc_system}
By $\mu$ we denote the standard measure on $S^1$. Let $g \colon S^1 \to [0,\infty)$ be a measurable function with $\int g \,\mathrm{d}\mu=1$ and let $\mathrm{d}\pi = g \mathrm{d}\mu$. Then $\pi$ is a probability measure and we have the following lower bound for the entropy of the pair $(\Ic_\alpha,\pi)$ for any $\alpha \in (0,1)$:
\[ H(\Ic_\alpha, \pi) \geq 
- \int g \log g \, \mathrm{d} \mu  - \log \alpha .\]
Moreover, this holds with equality if $g$ can be obtained as 
\[ g(x) = \frac{1}{\alpha} \nu\big( (x,x+\alpha) \big) 
\mbox{ for a.e.\ $x \in S^1$ for some Borel probability measure $\nu$;} \]
in particular, whenever there exists a measurable $\hat{g} \colon S^1 \to [0,\infty)$ such that 
\[ g(x) = \frac{1}{\alpha} \int_x^{x+\alpha} \hat{g} \,\mathrm{d}\mu 
\mbox{ for a.e.\ } x \in S^1 .\]
\end{proposition}
\begin{proof}
The function $b=\alpha g$ satisfies condition \eqref{eq:condition}:
\[ \int_x^{x+\alpha} \frac{1}{b} \,\mathrm{d}\pi = 
\int_x^{x+\alpha} \frac{1}{\alpha g} g \,\mathrm{d}\mu = 1 .\]
This gives us the lower bound 
\[ \fpi(b) = - \int \log g +  \log \alpha \,\mathrm{d}\pi = 
- \int g \log g \, \mathrm{d} \mu  - \log \alpha .\]

If the extra condition is satisfied, then the above $b$ is clearly in $\VPfull(\Ic_\alpha,\pi)$, and hence $\fpi(b)$ is an upper bound as well.
\end{proof}

\subsection{System of sets of small measure}
Our lower bound says that $H(\Ic,\pi) \geq \fpi(a)$ provided that the integral of $1/a$ over each set in $\Ic$ is at most $1$. In a reverse manner, let us start with a fixed measurable function $a \colon \Omega \to (0,1]$ and define $\Ic$ as the collection of all sets over which the integral of $1/a$ is at most $1$. 
\begin{example} 
Suppose that $(\Omega,\Ac,\pi)$ is a probability space. For a fixed measurable function $a \colon \Omega \to (0,1]$ consider the following measure
\[ \mu(A) = \int_A \frac{1}{a} \,\mathrm{d}\pi; \mbox{ in notation: } 
\mathrm{d}\mu = \frac{1}{a} \mathrm{d}\pi .\]
Then let 
\[ \Ic \defeq \big\{ J \in \Ac \, : \, \mu(J) \leq 1 \big\} 
= \bigg\{ J \subseteq \Omega \mbox{ measurable} \, : \, 
\int_J \frac{1}{a} \, \mathrm{d} \pi \leq 1 \bigg\}. \]
\end{example}
\begin{proposition}
We always have $H(\Ic, \pi) \geq \fpi(a)$. Moreover, if $\pi$ is atomless, we have equality:
\[ H(\Ic, \pi) = \fpi(a) .\]
Consequently, if 
\[ \mu(\Omega) = \int \frac{1}{a} \,\mathrm{d}\pi=\infty \mbox{ but } 
\fpi(a) = \int \log\left(\frac{1}{a}\right) \,\mathrm{d}\pi < \infty ,\]
then $\Ic$ has infinite chromatic number but finite entropy. 
\end{proposition}
\begin{proof}
By construction, $a$ satisfies the condition of the lower bound. The upper bound will follow easily from the next lemma. We omit the proof which is fairly standard. 
\begin{lemma}
Let $(\Omega,\Ac,\mu)$ be an atomless measure space. Suppose that $f \colon \Omega \to [0,\infty)$ is a measurable function such that its $L_1(\mu)$ and $L_\infty(\mu)$ norms satisfy that 
\[ \| f \|_{\infty} \leq \| f \|_{1} < \infty .\]
Then there exist measurable sets $J_k$ and nonnegative coefficients $q_k \geq 0$ such that
\[ \forall k\in\Nb \ \mu(J_k)=1; \sum_{k=1}^\infty q_k = \| f \|_{1} ;
\mbox{ and $\mu$-a.e.\ }
f = \sum_{k=1}^\infty q_k \ind_{J_k} .\]
\end{lemma}
Setting $\ds \mathrm{d}\mu=\frac{1}{a} \mathrm{d}\pi$ and $f=a$, we have $\|f\|_1=1$ and the lemma yields that $a \in \VPcnt(\Ic)$, showing that $\fpi(a)$ is an upper bound for the entropy.
\end{proof}

\subsection{System of independent events} \label{sec:indep}
In the following family of examples the chromatic number is infinite, while the entropy is typically finite.

It can be seen easily that for any system $\Ic$, one can choose a countable subsystem $\Ic'=\{J_1, J_2, \ldots\} \subseteq \Ic$ such that 
\[ \Hlow_\la(\Ic',\pi) = \Hlow_\la(\Ic,\pi) \mbox{ and } 
\Hupp_\la(\Ic',\pi) = \Hupp_\la(\Ic,\pi) \mbox{ for all } \la \in (0,1) .\]
So determining the entropy for countable systems is certainly a key problem. For $\Ic=\{J_1, J_2, \ldots\}$, the problem is described by the probabilities $\pi( J_{i_1} \cap \cdots \cap J_{i_\ell} )$, running through finite sequences of indices $i_1 < \cdots < i_\ell$. From this viewpoint, a very natural special case is when $J_1, J_2, \ldots$ are mutually independent (as random events in $\Omega$). In this case the only relevant parameters are 
\[ m_k \defeq \pi(J_k) .\]
\begin{example} \label{ex:indep}
Given $m_k \in (0,1)$, $k=1,2,\ldots$, let $J_1, J_2, \ldots$ be mutually independent events in a (suitable) probability space $(\Omega,\Ac,\pi)$ with $\pi(J_k)=m_k$, and consider the system $\Ic=\{J_1, J_2, \ldots\}$.

In particular, if $m_k=1/2$ for each $k$, then one can take $\Omega=[0,1)$ with the Lebesgue measure $\pi$, and $J_k$ consisting of $x \in [0,1)$ for which the $k$-th binary digit is $1$, that is, $J_1=[1/2,1)$; $J_2=[1/4,1/2) \cup [3/4,1)$; and so on. 
\end{example}
\begin{proposition} \label{prop:half}
If $m_k=1/2$ for each $k$, then 
\[ H(\Ic,\pi) = \log(2) .\]
In other words, it has the same entropy as the larger system $\Ic'$ consisting of all measurable sets $J \subseteq \Omega$ with $\pi(J) \leq 1/2$.
\end{proposition}
\begin{proof}
Since $\pi(J_k)=1/2$ for each $k \geq 1$, the constant function $b=1/2$ satisfies the condition of the lower bound, so we get $H(\Ic,\pi) \geq \log(2)$. 
Next we will find countable convex combinations $a \in \VPcnt(\Ic)$ with $\fpi(a)$ approaching $\log 2$. For the sake of simplicity, we will use the particular setting $\Omega=[0,1)$ and $J_k$ being the set of numbers whose $k$-th digit is $1$.

As a first step. we show the existence of $a_0 \in \VPcnt(\Ic)$ with $\fpi(a_0)<\infty$. Notice that for 
\[ a_0 \defeq \sum_{k=1}^\infty \frac{1}{2^k} \ind_{J_k} \in \VPcnt(\Ic) \]
we have $a_0(x)=x$ for all $x\in[0,1)$, and hence $\fpi(a_0)=-\int_0^1 \log x \,\mathrm{d}x = 1$.

For a given $n \in \Nb$ let 
\[ a_n \defeq \frac{1}{n} \sum_{k=1}^n \ind_{J_k} \in \VP(\Ic) .\]
Since $J_1, J_2, \ldots$ are mutually independent, it follows from the law of large numbers that $a_n$ converges to the constant $1/2$ function as $n \to \infty$. However, $a_n(x)=0$ on a positive measure set, and hence $\fpi(a_n) = \infty$ for all $n$. We can fix this by mixing $a_n$ with $a_0$: for $\eps>0$ let 
\[ a_{n,\eps} \defeq \eps a_0 + (1-\eps)a_n \in \VPcnt(\Ic) .\]
Furthermore, we define the following sets where the ``average'' $a_n$ is smaller than the typical value $1/2$:
\[ A_{n,\de} \defeq \big\{ x \, : \, a_n(x) < \frac{1}{2} - \de \big\} .\]
The weak law of large numbers implies that $\pi( A_{n,\de} ) \to 0$ as $n \to \infty$ for any fixed $\de>0$. For $x \notin A_{n,\de}$ we can use the bound $a_{n,\eps}(x) \geq (1/2 - \de)(1-\eps)$, while for $x \in A_{n,\de}$ we have $a_{n,\eps}(x) \geq \eps a_0(x)$. It follows that 
\[ \fpi( a_{n,\eps} ) \leq -\log\big( (1/2 - \de)(1-\eps) \big) - \int_{A_{n,\de}} \log(\eps a_0) .\]
Since $a_0(x)=x$, the second term can be bounded by 
\[ \int_0^{\pi(A_{n,\de})} - \log(\eps x) \,\mathrm{d}x = 
\pi(A_{n,\de}) \big( -\log \eps - \log \pi(A_{n,\de}) + 1 \big) .\]
Taking $n \to \infty$, and then $\delta, \eps \to 0$ we get that 
\[ H( \Ic , \pi ) \leq 
\inf_{n, \eps>0} \fpi( a_{n,\eps} ) \leq - \log(1/2) = \log 2 ,\]
which completes the proof.
\end{proof}
One could use a similar argument to settle the slightly more complicated case when all probabilities but one are the same.
\begin{proposition}
For $1> m_1 \geq m_\infty >0$ let $J_1, J_2, \ldots$ be independent events with 
\[ \pi(J_1) = m_1 \mbox{ and } 
\pi(J_k)=m_\infty \mbox{ for all } k \geq 2 .\]
Then we have the following formula for the entropy of $\Ic=\{J_1, J_2, \ldots\}$:
\begin{equation} \label{eq:allbutone}
H(\Ic,\pi) = -m_1 \log m_1 +  
(1-m_1)\big( \log(1-m_\infty) - \log m_\infty - \log(1-m_1) \big) .
\end{equation}
Furthermore, note that $H(\Ic,\pi) \leq H( \{J_2,J_3,\ldots\} , \pi ) = - \log m_\infty$.
\end{proposition}
\begin{proof}[Sketch of the proof]
Set the coefficients as
\[ q_1 \defeq \frac{m_1-m_\infty}{1-m_\infty} \mbox{ and }
q_\infty \defeq 1-q_1 = \frac{1-m_1}{1-m_\infty} \]
and consider the corresponding convex combination
\[ a \defeq q_1 \ind_{J_1} + q_\infty m_\infty = 
m_1 \ind_{J_1} + \frac{m_\infty(1-m_1)}{1-m_\infty}\ind_{\Omega \sm J_1} .\]
Using similar methods as in the proof of Proposition \ref{prop:half}, the entropy can be shown to be equal to $\fpi(a)$, which gives \eqref{eq:allbutone}.

When $m_1=m_\infty$, that is, when all events have the same probability $m_\infty$, then the formula gives  $-\log m_\infty$. This is actually an upper bound for \eqref{eq:allbutone} because the entropy of a subsystem is clearly greater than or equal to the entropy of the system, implying $H(\Ic,\pi) \leq H( \{J_2,J_3,\ldots\} , \pi ) = - \log m_\infty$.
\end{proof}
\begin{remark}
Example \ref{ex:indep} is interesting even for probabilities $m_k$ converging to $0$. A necessary condition for the entropy to be finite is 
\[ \pi\left( \bigcup_{k=1}^\infty J_k \right) = 1 
\Longleftrightarrow \prod_{k=1}^\infty (1-m_k) = 0 
\Longleftrightarrow \sum_{k=1}^\infty m_k = \infty .\]
\end{remark}
So an interesting choice of parameters would be, for example, $\ds m_k=\frac{1}{k+1}$. Considering 
\[ a \defeq \sum_{k=1}^\infty q_k \ind_{J_k} \mbox{ with coefficients }
q_k= \frac{1}{k(k+1)} ,\]
we get the upper bound 
\[ \fpi(a) \leq 
\sum_{k=1}^\infty (1-m_1)\cdots(1-m_{k-1}) m_k \log( 1/q_k ) = 
\sum_{k=1}^\infty \frac{\log k + \log(k+1)}{k(k+1)} < \infty .\]
It follows that the entropy is finite but we could not determine its precise value.

\subsection{Translated copies}
\begin{example} \label{ex:tr_copies}
Let $\Omega=[0,1)$ with the Lebesgue measure $\pi$. For a measurable $A \subset [0,1)$ with $\pi(A)=\al>0$ we consider the translated copies of $A$ modulo $1$:
\[ \Ic \defeq \big\{ [0,1) \cap \big( (A+t) \cup (A+t-1) \big) 
\, : \, t \in [0,1) \big\} .\]
The entropy of the system is $-\log \al$ because the constant $\al$ function lies in $\VPfull(\Ic,\pi)$ and it clearly satisfies the condition of the lower bound as well. Therefore:
\[ H(\Ic,\pi) = \Hlow(\Ic,\pi) = \Hfull(\Ic, \pi) = - \log \al .\]
However, the chromatic number is not necessarily finite. For example, suppose that $A$ is a Cantor set of positive measure (i.e., a nowhere dense perfect set with positive measure). Then a finite union of translated copies of $A$ is still closed and nowhere dense, therefore it cannot have full measure, proving that the chromatic number is not finite. In other words, any function $a \in \VP(\Ic)$ is $0$ on a set of positive measure, and hence $\Hfin(\Ic,\pi)=\infty$.
\end{example}
As for $\Hcnt(\Ic,\pi)$, it seems to be equal to $-\log \al$ for any $A$. In fact, the following may be true in general.
\begin{conjecture} \label{conj:full_vs_cnt}
For any system $\Ic$ we have 
\[ \Hfull(\Ic,\pi)=\Hcnt(\Ic,\pi) .\] 
Moreover, if $\chi(\Ic) < \infty$, then 
\[ \Hfull(\Ic,\pi)=\Hcnt(\Ic,\pi)=\Hfin(\Ic,\pi)<\infty .\]
\end{conjecture}

\bibliographystyle{plain}
\bibliography{refs}

\begin{thebibliography}{10}

\bibitem{aldous}
D.~J. Aldous.
\newblock Exchangeability and related topics.
\newblock In {\em École d'Été de Probabilités de Saint-Flour XIII, 1983},
  volume 1117, pages 1--198. Springer, 1985.

\bibitem{entropy_splitting}
Imre Csisz\'ar, J\'anos K\"orner, L\'aszl\'o Lov\'asz, Katalin Marton, and
  G\'abor Simonyi.
\newblock Entropy splitting for antiblocking corners and perfect graphs.
\newblock {\em Combinatorica}, 10(1):27--40, 1990.

\bibitem{hatami_norine_2013}
Hamed Hatami and Serguei Norine.
\newblock The entropy of random-free graphons and properties.
\newblock {\em Combinatorics, Probability and Computing}, 22(4):517–526,
  2013.

\bibitem{hladky_rocha}
Jan Hladký and Israel Rocha.
\newblock Independent sets, cliques, and colorings in graphons.
\newblock {\em European Journal of Combinatorics}, 88:103--108, 2020.
\newblock Selected papers of EuroComb17.

\bibitem{janson}
Svante Janson.
\newblock {\em Graphons, cut norm and distance, couplings and rearrangements},
  volume~4 of {\em New York Journal of Mathematics. NYJM Monographs}.
\newblock State University of New York, University at Albany, Albany, NY, 2013.

\bibitem{sorting}
J.~Kahn and J.H. Kim.
\newblock Entropy and sorting.
\newblock {\em Journal of Computer and System Sciences}, 51(3):390--399, 1995.

\bibitem{hashing}
J.~K{\"o}rner and K.~Marton.
\newblock New bounds for perfect hashing via information theory.
\newblock {\em European Journal of Combinatorics}, 9(6):523--530, 1988.

\bibitem{korner1973}
J{\'a}nos K{\"o}rner.
\newblock Coding of an information source having ambiguous alphabet and the
  entropy of graphs.
\newblock In {\em 6th Prague conference on information theory}, pages 411--425,
  1973.

\bibitem{leader}
Imre Leader.
\newblock The fractional chromatic number of infinite graphs.
\newblock {\em Journal of Graph Theory}, 20(4):411--417, 1995.

\bibitem{lovaszbook}
L.~Lov{\'a}sz.
\newblock {\em Large Networks and Graph Limits}.
\newblock American Mathematical Society colloquium publications. American
  Mathematical Society, 2012.

\bibitem{circuit}
I.~Newman, P.~Ragde, and A.~Wigderson.
\newblock Perfect hashing, graph entropy, and circuit complexity.
\newblock In {\em Proceedings Fifth Annual Structure in Complexity Theory
  Conference}, pages 91--99, 1990.

\bibitem{boolean}
Jaikumar Radhakrishnan.
\newblock {$\Sigma\Pi\Sigma$} threshold formulas.
\newblock {\em Combinatorica}, 14(3):345--374, 1994.

\bibitem{survey}
G{\'a}bor Simonyi.
\newblock Graph entropy: A survey.
\newblock In William Cook, L{\'a}szl{\'o} Lov{\'a}sz, and Paul Seymour,
  editors, {\em Combinatorial Optimization}, volume~20 of {\em DIMACS Series in
  Discrete Mathematics and Theoretical Computer Science}, pages 399--441, 1993.

\bibitem{survey2}
G{\'a}bor Simonyi.
\newblock Perfect graphs and graph entropy. {A}n updated survey.
\newblock In Jorge Ramirez-Alfonsin and Bruce Reed, editors, {\em Perfect
  Graphs}, pages 293--328. John Wiley and Sons, 2001.

\end{thebibliography}

\end{document}